\documentclass[11pt]{article}
\usepackage{amsfonts,latexsym,amsmath,amscd,geometry}
\date{\empty}

\usepackage{amsmath}
\usepackage{amssymb}
\usepackage{latexsym}
\newcommand \nc {\newcommand}
\nc \proof {\noindent {\em{Proof.\/ }}}
\nc \qed {\hfill $\Box$}

\newtheorem{Theorem}{Theorem}[section]
\newtheorem{Lemma}{Lemma}[section]

\newtheorem{definition}[Theorem]{Definition}

\newcommand{\beq}{\begin{eqnarray}}
\newcommand{\eeq}{\end{eqnarray}}
\newcommand{\beqno}{\begin{eqnarray*}}
\newcommand{\eeqno}{\end{eqnarray*}}
\newcommand{\be}{\begin{equation}}
\newcommand{\ee}{\end{equation}}

\def \d2{\Delta_{+}^2}
\begin{document}
\title{Well-posedness for the heat flow of harmonic maps and the liquid crystal flow with rough initial data}
\author{Changyou Wang
\thanks{Department of Mathematics, University of Kentucky
Lexington, KY 40506}}
\date{}
\maketitle
\begin{abstract}
We investigate the well-posedness of (i) the heat flow
of harmonic maps from $\mathbb R^n$ to a compact Riemannian manifold $N$ without boundary for initial data
in BMO;
and (ii) the hydrodynamic flow $(u,d)$ of nematic liquid crystals on $\mathbb R^n$
for initial data in ${\rm{BMO}}^{-1}\times{\rm{BMO}}$.
\end{abstract}

\setcounter{section}{0} \setcounter{equation}{0}
\section{Introduction}

For $k\ge 1$, let $N$ be a $k$-dimensional
compact Riemannian manifold without boundary, isometrically
embedded in some Euclidean space $\mathbb R^l$.
For $n\ge 1$, the equation of heat flow of harmonic maps 
from $\mathbb R^n$ to $N$ is given by:
\begin{eqnarray} \label{hfhm}
u_t-\Delta u&=&A(u)(\nabla u,\nabla u) \ \ {\rm{in}}\ \mathbb R^n\times (0,+\infty)\\
u\big|_{t=0}&=& u_0 \ \ \ {\rm{in}}\ \mathbb R^n \label{ic}
\end{eqnarray}
where $A(y): T_y N\times T_y N\to (T_y N)^\perp$ is the second fundamental form of $N\subset\mathbb R^l$
at $y\in N$, and $u_0:\mathbb R^n\to N$ is a given map. 

(\ref{hfhm})-(\ref{ic})  provides a very important approach to seek the existence of harmonic maps in various topological
classes. In their pioneering work \cite{ES} in 1960's, Eells-Sampson established that (i) for $u_0\in C^\infty(\mathbb R^n, N)$ 
there exists $0<T=T(\phi)\leq+\infty$ such that (\ref{hfhm})-(\ref{ic}) admits a unique smooth 
solution $u\in C^\infty(\mathbb R^n\times [0,T),N)$; and (ii) if, in additions, the sectional
curvature $K_N$ of $N$ is nonpositive, then $u\in C^\infty(\mathbb R^n\times \mathbb R_+, N)$ and 
\begin{equation}
\|u\|_{C^2(\mathbb R^n\times \mathbb R_+)}\le C(n, \|\phi\|_{C^2(\mathbb R^n)}). \label{prior_est} 
\end{equation}
Without the curvature assumption,
Hildebrandt-Kaul-Widman \cite{HKW} established, in late 70's, the existence of unique, global smooth solution to
(\ref{hfhm})-(\ref{ic}) under the assumption that the image of $u_0$ is contained in a geodesic ball $B_R$ in $N$
with radius $R<\frac{\pi}{2\sqrt{\max_{B_R}|K_N|}}$. In general, on the one hand, 
it is well-known via the works by Coron-Ghidaglia \cite{CG}, Chen-Ding \cite{CD}, and Chang-Ding-Ye \cite{CDY} that
the short time smooth solution to (\ref{hfhm})-(\ref{ic}) may develop finite time singularity;
on the other hand, Chen-Struwe \cite{CS} (see also Chen-Lin \cite{CL}) established the existence
of partially smooth, global weak solutions to (\ref{hfhm})-(\ref{ic}) for smooth initial data
$u_0$.

Although there have been many important works to (\ref{hfhm})-(\ref{ic}) (see for example
Lin-Wang \cite{LW} and references therein), it remains an interesting question 
the global (or local, resp.) well-posedness of (\ref{hfhm})-(\ref{ic}) for small 
(or large, resp.) rough initial data.  For initial data $u_0$ in the Sobolev space
$W^{1,p}(\mathbb R^n, N)$ ($1<p\leq +\infty$), Struwe \cite{Struwe} established, in
dimension $n=2$, the local well-posedness 
of (\ref{hfhm})-(\ref{ic}) in the space $L^2_t H^2_x$ for
$u_0\in W^{1,2}(\mathbb R^2, N)$, and the global well-posedness provided
$\|\nabla u_0\|_{L^2(\mathbb R^2)}$ is sufficiently small. 
For $n\ge 3$,  the well-posedness similar to that of \cite{Struwe} for 
$u_0\in W^{1,n}(\mathbb R^n,N)$ was not available in the literature 
previously, and the
readers can refer to Wang \cite{Wang} for some related earlier results. 

In a very interesting paper \cite{KL}, Koch-Lamm proved that (\ref{hfhm})-(\ref{ic})
is (i) locally uniquely solvable in $C^\infty(\mathbb R^n,N)$ provided $u_0$
is $L^\infty$-close to a uniformly continuous map; and (ii) globally uniquely
solvable in $C^\infty(\mathbb R^n,N)$ provided $u_0$ is $L^\infty$-close to a point.
The techniques employed by Koch-Lamm in \cite{KL} were originated from the
earlier work by Koch and Tataru \cite{KT} on the global 
well-posedness of the incompressible Navier-Stokes
equation for $u:\mathbb R^n\times\mathbb R_+\to\mathbb R^n$:
\begin{eqnarray} \label{nse}
 u_t+u\cdot\nabla u-\Delta u+\nabla P &=& 0   \ \ {\rm{in}}\ \mathbb R^n\times (0,+\infty) \\
\nabla \cdot u&=& 0  \ \ \ {\rm{in}} \ \mathbb R^n \label{incomp}\\
u\big|_{t=0}&=& u_0 \ \ \ {\rm{in}}\ \mathbb R^n \label{ic1}
\end{eqnarray}
for $u_0\in {\rm{BMO}}^{-1}(\mathbb R^n)$ with $\nabla\cdot u_0=0$ and small $\|u_0\|_{{\rm{BMO}}^{-1}}$.

Partially motivated by \cite{KL} and \cite{KT}, we address the well-posedness for both the heat
flow of harmonic maps and the hydrodynamic flow of nematic liquid crystals in this paper.

In order to state the results, we first recall the definitions of both the local and global
${\rm{BMO}}$  spaces.
\begin{definition} For $0<R\leq +\infty$, a function $f\in L^1_{\rm{loc}}(\mathbb R^n)$
is in ${\rm{BMO}}_R(\mathbb R^n)$ if the semi-norm
$$
\left[f\right]_{{\rm{BMO}}_R(\mathbb R^n)}
:=\sup_{x\in\mathbb R^n, 0<r\le R}\left\{ r^{-n}\int_{B_r(x)}|f(y)-f_{x,r}|\,dy\right\}
$$
is finite, where $f_{x,r}=\frac{1}{|B_r(x)|}\int_{B_r(x)} f(y)\,dy$ is the average of $f$ over $B_r(x)$. 
We say $f\in \overline{\rm{VMO}}(\mathbb R^n)$ if
$$\lim_{r\downarrow 0}\left[f\right]_{{\rm{BMO}}_r(\mathbb R^n)}=0.
$$
When $R=+\infty$, we simply write $({\rm{BMO}}(\mathbb R^n), [\cdot]_{{\rm BMO}(\mathbb R^n)})$ 
for $({\rm{BMO}}_{\infty}(\mathbb R^n), [\cdot]_{{\rm BMO}_\infty(\mathbb R^n)})$.
\end{definition}

Now recall the space ${\rm{BMO}}^{-1}$, introduced  by Koch-Tataru \cite{KT},  as follows.

\begin{definition} For $0<R\leq +\infty$, a function $f\in L^1_{\rm{loc}}(\mathbb R^n)$
is in ${\rm{BMO}}_R^{-1}(\mathbb R^n)$ if there exists $(f_1,\cdots, f_n)\in {\rm{BMO}}_R(\mathbb R^n)$
such that $f=\sum_{i=1}^n \frac{\partial f_i}{\partial x_i}$. Moreover, the norm
of $f$ is defined by
$$
\left\|f\right\|_{{\rm{BMO}}_R^{-1}(\mathbb R^n)}
:=\inf\left\{ \sum_{i=1}^n [f_i]_{{\rm{BMO}}_R(\mathbb R^n)}: \  f\equiv\sum_{i=1}^n \frac{\partial f_i}{\partial x_i}\right\}.
$$
We say $f\in (\overline{\rm{VMO}}(\mathbb R^n))^{-1}$ if
$$\lim_{r\downarrow 0}\left[f\right]_{{\rm{BMO}}_r^{-1}(\mathbb R^n)}=0.
$$
When $R=+\infty$, we simply write $({\rm{BMO}}^{-1}(\mathbb R^n), [\cdot]_{{\rm BMO}^{-1}(\mathbb R^n)})$ for $({\rm{BMO}}_{\infty}^{-1}(\mathbb R^n), [\cdot]_{{\rm BMO}^{-1}_\infty(\mathbb R^n)})$.
\end{definition}

We also introduce the functional space $X_T$ for $0<T\leq +\infty$ as follows.
$$
X_T:=\left\{f:\mathbb R^n\times [0,T]\to \mathbb R^l \ \big | \ \ 
|||u|||_{X_T}\equiv \sup_{0<t\le T}\|f(t)\|_{L^\infty(\mathbb R^n)}+\left\|f\right\|_{X_T}<+\infty\right\},
$$
where
$$\left\|f\right\|_{X_T}=\sup_{0<t\le T} \sqrt{t}\|\nabla f(t)\|_{L^\infty(\mathbb R^n)}
+\sup_{x\in\mathbb R^n, 0<R\le\sqrt{T}} \ (R^{-n}\int_{P_R(x,R^2)}|\nabla f|^2\,dxdt)^\frac12,$$
and $P_R(x,R^2)=B_R(x)\times [0,R^2]$ denotes the parabolic cylinder with center  $(x,R^2)$
and radius $R$. It is easy to see that $(X_T, |||\cdot|||_{X_T})$ is a Banach space.
When $T=+\infty$, we simply write $X$ for $X_{\infty}$, $\|\cdot\|_X$ for $\|\cdot\|_{X_\infty}$,
and $[\cdot]_{X}$ for $[\cdot]_{X_\infty}$ respectively. 

For the heat flow of harmonic maps, we prove 
\begin{Theorem} {\rm{(local well-posedness)}}\label{lwp}
There exists $\epsilon_0>0$ such that for any $R>0$ if $\|u_0\|_{{\rm{BMO}}_R(\mathbb R^n)}\le\epsilon_0$,
then (\ref{hfhm})-(\ref{ic})  has a unique solution
$u\in X_{R^2}$ with small $\|u\|_{X_{R^2}}$.
In particular, if $u_0\in \overline{\rm{VMO}}(\mathbb R^n)$, then there exists $T_0>0$ such that
(\ref{hfhm})-(\ref{ic}) admits a unique solution $u\in X_{T_0}$ with small $\|u\|_{X_{T_0}}$.
\end{Theorem}

As a corollary, we have 
\begin{Theorem}{\rm{(global well-posedness)}}\label{gwp} There exist $\epsilon_0>0$ and
$C_0>0$ such that
if $[u_0]_{\rm{BMO}(\mathbb R^n)}\le\epsilon_0$, then there exists a unique
global solution $u\in X$ to
(\ref{hfhm})-(\ref{ic}) such that $\|u\|_X\le C_0\epsilon_0.$
\end{Theorem}

Since $W^{1,n}(\mathbb R^n)\subset\overline{{\rm{VMO}}}(\mathbb R^n)$, it follows
from Theorem 1.3 that for any initial data $u_0\in W^{1,n}(\mathbb R^n)$, (\ref{hfhm})-(\ref{ic})
admits a short time unique solution $u\in X_{T_0}$ for some $T_0>0$. 
Theorem 1.4 implies that such a unique solution $u$ is a unique global solution in $X$ provided 
$\|\nabla u_0\|_{L^n(\mathbb R^n)}$ is sufficiently small.

Now we turn to the discussion on the well-posedness for the hydrodynamic flow of nematic liquid crystals in the entire space. 

The following equation modeling the hydrodynamic flow of namatic liquid crystal materials
has been proposed and investigated by  Lin-Liu \cite{LL1, LL2} in 1990's. 
\begin{eqnarray} \label{llf1}
u_t +u\cdot \nabla u-\Delta u+\nabla P&=&-\nabla\cdot(\nabla d\otimes\nabla d) 
\ \ \ \  \ \ \ \ {\rm{in}} \ \mathbb R^n\times (0,+\infty)\\
\nabla\cdot u &=& 0
\ \ \ \ \ \ \  \ \ \ \ \ \ \ \ \ \ \ \ \ \ \  \ \ \ \  \ {\rm{in}} \ \mathbb R^n\times (0,+\infty) \label{llf2}\\
d_t +u·\cdot \nabla d&=& \Delta d+|\nabla d|^2d 
\ \ \ \ \ \ \ \ \ \ \ \ \  {\rm{in}} \ \mathbb R^n\times (0,+\infty),
\label{llf3}
\end{eqnarray}
where $u(\cdot, t) : \mathbb R^n \to\mathbb R^n $
represents the velocity field of the flow, 
$d(\cdot, t) : \mathbb R^n \to S^2$, the unit sphere in $\mathbb R^3$, is a unit-vector field that represents the macroscopic molecular orientation of the nematic 
liquid crystal material, and $P(\cdot, t) : \mathbb R^n \to \mathbb R$
represents the pressure function.
$\nabla\cdot$ denotes the divergence operator, and 
$\nabla d \otimes \nabla d$ denotes the  $n \times n$
matrix whose $(i, j)$-the entry is given by $\nabla_i d \cdot \nabla_j d$
for $1 \leq i, j \leq n$.

The above system is a simplified version of the Ericksen-Leslie model, which 
reduces to the Ossen-Frank model in the static case, for the hydrodynamics of nematic
liquid crystal materials developed during the period of 1958 through 1968 
(see \cite{Ericksen, Gennes, Leslie}).  It is a macroscopic continuum description of the time evolution of the materials under the influence of both the flow field $u(x,t)$, and the macroscopic description of the microscopic orientation configurations $d(x,t)$ of rod-like liquid crystals. Roughly speaking, the system (\ref{llf1})-(\ref{llf3}) is a coupling between the incompressible Navier-Stokes equation and the transported heat flow of harmonic maps into $S^2$. 

When considering the initial and boundary value problem of (\ref{llf1})-(\ref{llf3}) on bounded domains $\Omega\subset\mathbb R^2$:
\begin{equation}
(u,d)\big|_{\Omega\times\{0\}}=(u_0,d_0), \ \ (u,d)\big|_{\partial \Omega\times (0,+\infty)}
=(0, d_0), \label{ibv}
\end{equation}
where $u_0:\Omega\to \mathbb R^2$ is a given divergence free vector field and
$d_0:\Omega\to S^2$ is a given unit-vector field. 
In a very recent paper, Lin-Lin-Wang \cite{LLW} proved, among other results, that for any
$(u_0,d_0)\in L^2(\Omega,\mathbb R^2)\times H^1(\Omega,S^2)$ with $\nabla\cdot u_0=0$, there is a global Leray-Hopf type weak solution $(u,d)$ to (\ref{llf1})-(\ref{llf3}) and
(\ref{ibv}) that is smooth away from at most finitely many singular times. 

In this paper, we want to address both local and global well-posedness issues on
the Cauchy problem of (\ref{llf1})-(\ref{llf3}) on $\mathbb R^n$ with rough initial data.

For this, we need to introduce another functional space in order to handle
the velocity field $u$. For $0<T\le +\infty$, let $Z_T$ be the space consisting
of functions $f:\mathbb R^n\times [0,T]$ such that
$$\|f\|_{Z_T}
:=\sup_{0<t\le T}\sqrt{t}\|f(t)\|_{L^\infty(\mathbb R^n)}
+\sup_{x\in\mathbb R^n, 0<r\le \sqrt{T}} (r^{-n}\int_{P_r(x,r^2)}|f|^2)^\frac12<+\infty.$$
When $T=+\infty$, we simply write $Z$ for $Z_\infty$, and $\|\cdot\|_Z$ for $\|\cdot\|_{Z_\infty}$. 

It turns out that, by combining the techniques of Koch-Tataru \cite{KT} and
Theorem 1.3 on the heat flow of harmonic maps, we are able to prove
the following theorems.
\begin{Theorem} There exists $\epsilon_0>0$ such that for any $R>0$ 
if $u_0\in {\rm{BMO}}_R^{-1}(\mathbb R^n, \mathbb R^n)$, with $\nabla\cdot u_0=0$,
and $d_0\in {\rm{BMO}}_R(\mathbb R^n,S^2)$ satisfies 
\beq{}\label{small}
\|u_0\|_{{\rm{BMO}}_R^{-1}(\mathbb R^n)}+\left[d_0\right]_{{\rm{BMO}}_R(\mathbb R^n)}
\le\epsilon_0,
\eeq
then there exists a unique solution $(u,d)\in Z_{R^2}\times X_{R^2}$ wth small 
$(\|u\|_{Z_{R^2}}+\|d\|_{X_{R^2}})$
to (\ref{llf1})-(\ref{llf3}) and 
\beq{} 
(u,d)\big|_{t=0}=(u_0,d_0)  \ \ {\rm{on}}\ \mathbb R^n. \label{llf4}
\eeq
In particular, if $(u_0,d_0)\in (\overline {\rm{VMO}}(\mathbb R^n))^{-1}
\times (\overline {\rm{VMO}}(\mathbb R^n))$, then 
there exists $T_0>0$ such that $(\ref{llf1})-(\ref{llf3})$ and  $(\ref{llf4})$ admits a unique solution
$(u,d)\in X_{T_0}$ with small $(\|u\|_{Z_{T_0}}+\|d\|_{X_{T_0}})$.
\end{Theorem}

As a corollary, we have 
\begin{Theorem} There exist $\epsilon_0>0$ and $C_0>0$ such that
if $u_0\in {\rm{BMO}}^{-1}(\mathbb R^n, \mathbb R^n)$, with $\nabla\cdot u_0=0$,
and $d_0\in {\rm{BMO}}(\mathbb R^n,S^2)$ satisfies 
\beq{}\label{small}
\|u_0\|_{{\rm{BMO}}^{-1}(\mathbb R^n)}+\left[d_0\right]_{{\rm{BMO}}(\mathbb R^n)}
\le\epsilon_0,
\eeq
then there exists a unique global solution $(u,d)\in Z\times X$ 
to (\ref{llf1})-(\ref{llf3}) and (\ref{llf4}) wth $(\|u\|_{Z}+\|d\|_{X})\leq C_0\epsilon_0$.
\end{Theorem}

We also remark that Theorem 1.5 implies that (\ref{llf1})-(\ref{llf3}) and (\ref{llf4})
is locally well-posed in $X_{T}$ for any initial data $(u_0,d_0)\in L^n(\mathbb R^n,
\mathbb R^n)\times W^{1,n}(\mathbb R^n,S^2)$, and is globally well-posed in $X$
provided $(\|u_0\|_{L^n(\mathbb R^n)}+\|\nabla d_0\|_{L^n(\mathbb R^n)})$
is sufficiently small. 

The remaining of this paper is organized as follows. In section 2, we establish
some basic estimates on the caloric extension of BMO functions. In section 3,
we prove Theorem 1.3 and 1.4. In section 4, we prove Theorem 1.5 and 1.6.

\setcounter{section}{1} \setcounter{equation}{0}
\section{Preliminary results}

In this section, we first review Carleson's well-known theorem on the characterization of
a BMO function in terms of its caloric extension, see Stein \cite{Stein} Page 159, Theorem 3.  Then we show a crucial estimate of the distance between the caloric extension of $u_0$ and the manifold $N$.
  
Let $G(x,t)$ be the fundamental solution of the heat equation in $\mathbb R^n\times \mathbb R_+$:
\beq{}\label{hk}
G(x,t)=\frac{1}{(4\pi t)^{\frac{n}2}} e^{-\frac{|x|^2}{4t}}, \ x\in\mathbb R^n, \ t>0.
\eeq
Let $\tilde{u}_0:\mathbb R^n\times \mathbb R_+\to\mathbb R^l$ be the caloric extension of $u_0$:
\beq{}\label{carloric}
\tilde{u}_0(x,t)=\int_{\mathbb R^n}G(x-y,t)u_0(y)\,dy.
\eeq
Carleson's characterization of the BMO space asserts that
$u_0\in {\rm{BMO}}(\mathbb R^n)$ iff $|\nabla\tilde u_0|^2\,dxdt$ is a Carleson measure on
$\mathbb R^n\times\mathbb R_+$, i.e.
$$\sup_{x\in\mathbb R^n, r>0}\ r^{-n}\int_{P_r(x,r^2)}|\nabla \tilde u_0|^2\,dxdt<+\infty,$$
and one has the equivalence of the norms:
\beq{}\label{bmo}
\left[u_0\right]_{\rm{BMO}(\mathbb R^n)}
\approx\sup_{x\in\mathbb R^n, r>0}\ (r^{-n}\int_{P_r(x,r^2)}|\nabla \tilde u_0|^2\,dxdt)^\frac12.
\eeq
If $u_0\in {\rm{BMO}}_R(\mathbb R^n)$ for some $0<R<+\infty$, 
then the same characterization as above gives
\beq{}\label{bmor}
\left[u_0\right]_{\rm{BMO}_R(\mathbb R^n)}
\approx\sup_{x\in\mathbb R^n, 0<r\leq R}\ (r^{-n}\int_{P_r(x,r^2)}|\nabla \tilde u_0|^2\,dxdt)^\frac12.
\eeq

Since $\tilde u_0$ solves the heat equation on $\mathbb R^n\times \mathbb R_+$, the standard gradient estimate implies that for any $t>0$,
\beq{}\label{grad_est}
\sqrt{t}\|\nabla \tilde u_0(t)\|_{L^\infty(\mathbb R^n)}
\lesssim \sup_{x\in\mathbb R^n}\ (t^{-\frac{n}2}\int_{P_{\sqrt t}(x,t)}|\nabla \tilde u_0|^2\,dyd\tau)^\frac12.
\eeq
In particular, we have that (i) if $u_0\in {\rm{BMO}}(\mathbb R^n)$, then
\beq{}\label{grad_est1}
\sup_{t>0} \sqrt{t}\|\nabla \tilde u_0\|_{L^\infty(\mathbb R^n)}
\lesssim \left[u_0\right]_{{\rm{BMO}}(\mathbb R^n)},
\eeq
and (ii) if $u_0\in {\rm{BMO}}_R(\mathbb R^n)$ for some $R>0$, then
\beq{}\label{grad_est2}
\sup_{0<t\leq R^2} \sqrt{t}\|\nabla \tilde u_0\|_{L^\infty(\mathbb R^n)}
\lesssim \left[u_0\right]_{{\rm{BMO}}_R(\mathbb R^n)}.
\eeq

Now we need to estimate the distance of $\tilde{u}_0$ to the manifold
$N$ in terms of the BMO norm of $u_0$, which plays an important role
in the proof of Theorems. More precisely, we have
\begin{Lemma} \label{distance} For any $\delta>0$, there exists $K=K(\delta,N)>0$ such that
if $u_0\in {\rm{BMO}}_R(\mathbb R^n)$ for some $0<R\leq+\infty$, then
\beq{}\label{distance1}
{\rm{dist}}(\tilde{u}_0(x,t),N)\le K \left[u_0\right]_{{\rm{BMO}}_R(\mathbb R^n)}+\delta,
\ \ \forall x\in\mathbb R^n, \ 0\leq t\le \frac{R^2}{K^2}.
\eeq
In particular, if $u_0\in {\rm{BMO}}(\mathbb R^n)$ then
\beq{}\label{distance2}
{\rm{dist}}(\tilde{u}_0(x,t),N)\le K \left[u_0\right]_{{\rm{BMO}}(\mathbb R^n)}+\delta,
\ \ \forall (x,t)\in\mathbb R^n\times \mathbb R_+.
\eeq
\end{Lemma}
\noindent{\it Proof}. Since (\ref{distance2}) follows directly from 
(\ref{distance1}) with $R=+\infty$, it suffices to prove (\ref{distance1}).
For any $x\in\mathbb R^n$, $t>0$, and $K>0$, denote
$$c_{x,t}^K=\frac{1}{|B_K(0)|}\int_{B_K(0)} u_0(x-\sqrt{t} z)\,dz.$$
Since
$$\tilde{u}_0(x,t)=\int_{\mathbb R^n} \frac{1}{(4\pi)^{\frac{n}2}} e^{-\frac{|y|^2}4} u_0(x-\sqrt{t} y)\,dy,$$
we have
\begin{eqnarray}\label{distance2}
\left|\tilde{u}_0(x,t)-c_{x,t}^K\right|
&\leq& \int_{\mathbb R^n} \frac{1}{(4\pi)^{\frac{n}2}} e^{-\frac{|y|^2}4} \left|u_0(x-\sqrt{t} y)-c_{x,t}^K\right|\,dy
\nonumber\\
&\leq&\left\{\int_{B_K(0)}+\int_{\mathbb R^n\setminus B_K(0)} \right\}
\frac{1}{(4\pi)^{\frac{n}2}} e^{-\frac{|y|^2}4} \left|u_0(x-\sqrt{t} y)-c_{x,t}^K\right|\,dy\nonumber\\
&\leq& \int_{B_K(0)}\left|u_0(x-\sqrt{t} y)-c_{x,t}^K\right|\,dy
+2\|u_0\|_{L^\infty(\mathbb R^n)}\int_{\mathbb R^n\setminus B_K(0)}e^{-\frac{|y|^2}4}\,dy\nonumber\\
&\leq& K^n \left[u_0\right]_{{\rm{BMO}}_{K\sqrt{t}}(\mathbb R^n)}+C_N\int_K^\infty e^{-\frac{r^2}4} r^{n-1}\,dr\nonumber\\
&\leq& \delta+ K^n \left[u_0\right]_{{\rm{BMO}}_{K\sqrt{t}}(\mathbb R^n)}
\end{eqnarray}
provided we choose a sufficiently large $K=K(\delta,N)>0$ so that 
$$C_N\int_K^\infty e^{-\frac{r^2}4} r^{n-1}\,dr\le\delta.$$
On the other hand, since $u_0(\mathbb R^n)\subset N$, we have
$$
{\rm{dist}}(c_{x,t}^K, N)\le \left|c_{x,t}^K-u_0(x-\sqrt{t} y)\right|, \ \forall y\in B_K(0)$$
and hence
\begin{equation}\label{distance3}
{\rm{dist}}(c_{x,t}^K, N)\le \frac{1}{|B_K(0)|}\int_{B_K(0)}|c_{x,t}^K-u_0(x-\sqrt{t} y)|\,dy
\le \left[u_0\right]_{{\rm{BMO}}_{K\sqrt{t}}(\mathbb R^n)}.
\end{equation}
Putting (\ref{distance2}) and (\ref{distance3}) together yields that (\ref{distance1}) 
holds for $t\le \frac{R^2}{K^2}$.
This completes the proof. \qed

\setcounter{section}{2} \setcounter{equation}{0}
\section{Proof of Theorem 1.3 and 1.4}

This section is devoted to the proof of Theorem 1.3 and 1.4. The idea is to choose
a suitable ball in $X$ such that the operator $T$ determined by the Duhamel formula 
has a fixed point in the ball.  

For $0<T\le +\infty$, besides the space $X_T$ introduced in section 1, we also
need to introduce $Y_T$ as follows.
$Y_T$ is the space consisting of all functions $f:\mathbb R^n\times [0,T]\to\mathbb R$
such that
$$\|f\|_{Y_T}\equiv \sup_{0<t\le T}\ t\|f(t)\|_{L^\infty(\mathbb R^n)}
+\sup_{x\in\mathbb R^n, 0<R\le\sqrt{T}} \ R^{-n}\int_{P_R(x,R^2)}|f|\,dxdt<+\infty.
$$
It is also easy to see $(Y_T, \|\cdot\|_{Y_T})$ is a Banach space. When $T=+\infty$,
we simply write $Y$ for  $Y_\infty$, and $\|\cdot\|_Y$ for $\|\cdot\|_{Y_\infty}$.

For $f\in Y_T$, define
\beq{}\label{resp}
\mathbb {S}f(x,t)=\int_0^t \int_{\mathbb R^n} G(x-y,t-s)f(y,s)\,dyds,\ (x,t)\in\mathbb R^n\times \mathbb R_+.
\eeq

It is well-known that $u:\mathbb R^n\times \mathbb R_+\to N$ solves (\ref{hfhm})-(\ref{ic})
iff
\beq{}\label{duhamel}
u(x,t)=\tilde u_0(x,t)+\mathbb {S}(A(u)(\nabla u,\nabla u))(x,t).
\eeq

The following Lemma plays the critical role in the proof.

\begin{Lemma} For $0<T\leq+\infty$, if $f\in Y_T$, then $\mathbb {S}f\in X_T$. Moreover,
\beq{}\label{YtoX}
|||\mathbb{S}f|||_{X_T}\le C\|f\|_{Y_T}
\eeq
for some $C=C(n)>0$.
\end{Lemma}
\noindent{\it Proof}.  By suitable scalings, we may assume $T\geq 1$.
Since the norms are invariant under both scaling and translation, it suffices to
show
\begin{equation}\label{yx_est}
|\mathbb{S}f(0,1)|+|\nabla (\mathbb{S}f)(0,1)|
+\left(\int_{P_1(0,1)}|\nabla (\mathbb{S}f)|^2\right)^\frac12\le C \|f\|_{Y_1}.
\end{equation}
Set $W=\mathbb{S}f$.  Then
\begin{eqnarray*}
W(0,1)&=&\int_0^1\int_{\mathbb R^n}G(y, 1-s) f(y,s)\,dyds\\
&=&\{\int_{\frac12}^1\int_{\mathbb R^n}+\int_0^\frac12\int_{B_2}+\int_0^\frac12\int_{\mathbb R^n\setminus B_2}\}
G(y, 1-s) f(y,s)\,dyds\\
&=&I_1+I_2+I_3.
\end{eqnarray*}
It is easy to see
$$|I_1|\leq \left(\sup_{\frac12\le s\le 1}\|f(s)\|_{L^\infty(\mathbb R^n)}\right)
\left(\int_{\frac12}^1\int_{\mathbb R^n}G(y,1-s)\,dyds\right)
\leq C\|f\|_{Y_1},$$
\begin{eqnarray*}
|I_2|&\leq& (\sup_{0\le s\le\frac12}\|G(\cdot,1-s)\|_{L^\infty(\mathbb R^n)})
(\int_{B_2\times [0,\frac12]}|f(y,s)|\,dyds)\\
&\leq& C\int_{B_2\times [0,\frac12]}|f(y,s)|\,dyds\leq C\|f\|_{Y_1},
\end{eqnarray*}
and
\begin{eqnarray*}
|I_3|&\leq&\int_0^\frac12 \int_{\mathbb R^n\setminus B_2}G(y,1-s)|f(y,s)|\,dyds\\
&\leq& C\int_0^\frac12\int_{\mathbb R^n\setminus B_2} e^{-\frac{|y|^2}{2}}|f(y,s)|\,dyds\\
&\leq& C \left(\sum_{k=2}^\infty k^{n-1} e^{-\frac{k^2}{2}}\right)
\cdot \left(\sup_{y\in \mathbb R^n} \int_{P_1(y,1)}|f(y,s)|\,dyds\right)\\
&\leq& C \|f\|_{Y_1}.
\end{eqnarray*}
Putting these three inequalities together implies $|W(0,1)|\le C\|f\|_{Y_1}$. The estimate of $|\nabla W(0,1)|$ can be
done similarly. In fact, denote
$$H(x,t)=\nabla_x G(x,t)=-\frac{x}{2t} G(x,t).$$
Then
$$\int_0^{\frac12}\int_{\mathbb R^n}|H(x,t)|\le C,
\ \sup_{x\in\mathbb R^n, \frac12\le t\le 1} |H(x,t)|\le C.$$
Since
$$\nabla W(0,1)=\int_0^1\int_{\mathbb R^n} H(-y, 1-s) f(y,s)\,dyds,
$$
we have
\begin{eqnarray*}
|\nabla W(0,1)|&\leq&  \int_0^1\int_{\mathbb R^n} |H|(-y, 1-s) |f(y,s)|\,dyds\\
&=&\{\int_{\frac12}^1\int_{\mathbb R^n}+\int_0^\frac12\int_{B_2}+\int_0^\frac12\int_{\mathbb R^n\setminus B_2}\}
|H(-y, 1-s)| |f(y,s)|\,dyds\\
&=&I_4+I_5+I_6.
\end{eqnarray*}
It is readily seen that
$$|I_4|\leq C(\int_0^\frac12\int_{\mathbb R^n}|H(x,t)|)\cdot(\sup_{\frac12\le s\le 1}\|f(\cdot,s)\|_{L^\infty(\mathbb R^n)})
\le C\|f\|_{Y_1},$$
$$|I_5|\leq C (\sup_{x\in\mathbb R^n, \frac12\le t\le 1} |H(x,t)|)(\int_{B_2\times [0,1]} |f(y,s)|\,dyds)
\le C\|f\|_{Y_1},$$
and
\begin{eqnarray*}
|I_6|&\leq& C\int_0^\frac12\int_{\mathbb R^n\setminus B_2}|y|e^{-\frac{|y|^2}{2}}|f(y,s)|\\
&\leq& C \left(\sum_{k=2}^\infty k^n e^{-\frac{k^2}{2}}\right)\cdot\left(\sup_{y\in\mathbb R^n}\int_{P_1(y,1)}|f(y,s)|\,dyds\right)\\
&\leq& C\|f\|_{Y_1}.
\end{eqnarray*}
Putting these estimates together yields $|\nabla W(0,1)|\le C\|f\|_{Y_1}.$

The estimate of $\|\nabla W\|_{L^2(P_1(0.1))}$ follows from the energy inequality as follows.
Since $W$ satisfies
$$W_t-\Delta W=f \ {\rm{in}}\ \mathbb R^n\times [0,1]; \ \ W\big|_{t=0}=0.$$
Let $\eta\in C_0^1(B_2)$ be a cut-off function of $B_1$. Multiplying the equation of $W$ by $\eta^2 W$ and
integrating over $\mathbb R^n\times [0,1]$, we obtain
\begin{eqnarray*}
\int_{P_1(0,1)}|\nabla W|^2 &\leq& C\int_{B_2\times [0,1]} (|W|^2+|W||f|)\\
&\leq& C\left(\|W\|_{L^\infty(B_2\times [0,1])}^2+ \|W\|_{L^\infty(B_2\times [0,1])}\|f\|_{L^1(B_2\times [0,1])}\right)\\
&\leq& C\|f\|_{Y_1}^2,
\end{eqnarray*}
where we have used in the last step the inequality, which was proved in the previous step,
$$\|W\|_{L^\infty(B_2\times [0,1])}\le C\|f\|_{Y_1}.$$
This completes the proof. \qed

In order to construct the solution to (\ref{hfhm}) in the space $X_{R^2}$,
we need to extend the second fundamental form $A(\cdot)(\cdot,\cdot)$ from $N$
to $\mathbb R^l$, still denoted as $A$. For this, recall that there exists $\delta_N>0$ such that
the nearest point projection map $\Pi:N_{\delta_N}=\{y\in\mathbb R^l:
{\rm{dist}}(y,N)\leq \delta_N\}\to N$ is smooth. Let $\widetilde{\Pi}\in C^\infty(\mathbb R^l,\mathbb R^l)$
be a smooth extension of $\Pi$, i.e. $\widetilde{\Pi}\equiv \Pi $ in $N_{\delta_N}$. Define
$$A(y)(V,W)=-D^2{\widetilde\Pi}(y)(V,W), \ \forall y\in \mathbb R^l, \ V, W\in T_y \mathbb R^l.$$

Now we define the mapping operator $\mathbf T$ on $X_{R^2}$ by letting
\begin{equation}\mathbf Tu(x,t)=\tilde u_0+\mathbb{S}(A(u)(\nabla u,\nabla u))(x,t), 
\ x\in\mathbb R^n, \  0<t\le R^2, \ \ u\in X_{R^2}. \label{t}
\end{equation}
If $u_0\in {\rm{BMO}}_R(\mathbb R^n)$,  then (2.4), (2.7), and the maximum principle
of the heat equation imply that $\tilde u_0\in X_{R^2}$ and
\beq{} \label{carloric_est}
\|\tilde {u}_0\|_{X_{R^2}}\lesssim \left[u_0\right]_{{\rm{BMO}}_R(\mathbb R^n)}.
\eeq

For $\epsilon>0$, let
$$\mathbf B_\epsilon(\tilde {u}_0)
:=\left\{u\in X: \big|\  \ |||u-\tilde{ u}_0|||_{X_{R^2}}\le\epsilon\right\}$$
be the ball in $X_{R^2}$ with center $\tilde{u}_0$ and radius $\epsilon$.
By the triangle inequality, we have
\beq{}\label{triangle}
\| u\|_{X_{R^2}}\le \|\tilde{u}_0\|_{X_{R^2}}
+\left\|u-\tilde{u}_0\right\|_{X_{R^2}}
\le \epsilon+\left\|\tilde{u}_0\right\|_{X_{R^2}}
\le\epsilon+\left[u_0\right]_{{\rm{BMO}}_R(\mathbb R^n)},
\ \forall u\in \mathbf B_\epsilon(\tilde{u}_0).
\eeq
In particular, we have
\begin{Lemma}. For $0<R\le +\infty$,  if $u_0:\mathbb R^n\to N$
satisfies $[u_0]_{{\rm{BMO}}_R(\mathbb R^n)}\le\epsilon$,
then
\beq{}\label{epsilon}
\|u\|_{L^\infty(\mathbb R^n\times [0,R^2])}\le C, \ \ 
\left\|u\right\|_{X_{R^2}}\le C\epsilon, \ \forall u\in \mathbf B_\epsilon(\tilde{u}_0)
\eeq
for some $C=C(n)>0$.
\end{Lemma}

Now we are ready to prove Theorem 1.3.  First we need the following two Lemmas.
\begin{Lemma} There exists $\epsilon_1>0$ 
such that if for $R>0$, $[u_0]_{{\rm{BMO}}_R(\mathbb R^n)}\le\epsilon_1$
then $\mathbf T$ maps $\mathbf B_{\epsilon_1}(\tilde u_0)$ to
$\mathbf B_{\epsilon_1}(\tilde u_0)$.
\end{Lemma}
\noindent{\it Proof}.
It follows from the formula (\ref{t}) that $\mathbf T(u)-\tilde u_0=
\mathbb{S}(A(u)(\nabla u,\nabla u))$
for $u\in \mathbf B_{\epsilon_1}(\tilde u_0)$. Hence Lemma 3.1 and Lemma 2.1 imply
\begin{eqnarray*}
&&|||\mathbf T(u)-\tilde{u}_0|||_{X_{R^2}}\\
&\leq& C\|A(u)(\nabla u,\nabla u)\|_{Y_{R^2}}\\
&=& C\left[\sup_{0<t\leq R^2} \ t\|A(u)(\nabla u,\nabla u)(t)\|_{L^\infty(\mathbb R^n)}
 +\sup_{x\in\mathbb R^n, 0<r\leq R} \ r^{-n}\int_{P_r(x,r^2)}|A(u)(\nabla u,\nabla u)|\right]\\
&\lesssim& \left(\sup_{0<t\leq R^2} \ \sqrt{t}\|\nabla u\|_{L^\infty(\mathbb R^n)}
 +\sup_{x\in\mathbb R^n, 0<r\leq R} \ (r^{-n}\int_{P_r(x,r^2)}|\nabla u|^2)^\frac12\right)^2\\
&\lesssim& \|u\|_{X_{R^2}}^2\leq C\epsilon_1^2\le\epsilon_1,
\end{eqnarray*}
provided  $\epsilon_1>0$ is chosen to be sufficiently small.  This completes the proof. \qed

\begin{Lemma} There exist $0<\epsilon_2\le\epsilon_1$ and $\theta_0\in (0,1)$
 such that if for $R>0$ 
$[u_0]_{{\rm{BMO}}_R(\mathbb R^n)}\le\epsilon_2$
then $\mathbf T:\mathbf B_{\epsilon_2}(\tilde u_0)\to \mathbf B_{\epsilon_2}(\tilde u_0)$
is a $\theta_0$-contraction map, i.e.
$$|||\mathbf T(u)-\mathbf T(v)|||_{X_{R^2}}\le\theta_0|||u-v|||_{X_{R^2}},
\ \forall u,v\in \mathbf B_{\epsilon_2}(\tilde {u}_0).
$$
\end{Lemma}
\noindent{\it Proof}. For $u,v\in \mathbf B_{\epsilon_2}(\tilde{u}_0)$, we have
\begin{eqnarray*}
 |\mathbf T u-\mathbf Tv|&=& |\mathbb{S}(A(u)(\nabla u,\nabla u)-A(v)(\nabla v,\nabla v))|\\
&\lesssim & |\mathbb{S}(A(u)(\nabla u,\nabla u)-A(u)(\nabla v,\nabla v))|
+|\mathbb{S}(A(u)(\nabla v,\nabla v)-A(v)(\nabla v,\nabla v))|\\
&\lesssim &\mathbb{S}((|\nabla u|+|\nabla v|)|\nabla(u-v)|)
+\mathbb{S}(|\nabla v|^2|u-v|).
\end{eqnarray*}
Hence, by Lemma 3.1, we obtain
$$
|||\mathbf T u-\mathbf Tv|||_{X_{R^2}}
\lesssim \|(|\nabla u|+|\nabla v|)|\nabla(u-v)|\|_{Y_{R^2}} 
+\||\nabla v|^2|u-v|\|_{Y_{R^2}}=I+II
$$
$I$ and $II$ can be estimated as follows.
\begin{eqnarray*}
&&I\\
&=&\sup_{0<t\le R^2} \ t\|(|\nabla u|+|\nabla v|)|\nabla(u-v)(t)|\|_{L^\infty(\mathbb R^n)}\\
&& +\sup_{x\in\mathbb R^n, 0<r\leq R} 
\ r^{-n}\int_{P_r(x,r^2)}(|\nabla u|+|\nabla v|)|\nabla(u-v)|\\
&\leq& \sup_{0<t\leq R^2} \ \sqrt{t}(\|\nabla u(t)\|_{L^\infty(\mathbb R^n)}+\|\nabla v(t)\|_{L^\infty(\mathbb R^n)})
\sup_{0<t\leq R^2} \ \sqrt{t}\|\nabla (u-v)(t)\|_{L^\infty(\mathbb R^n)}\\
&+&\sup_{x\in\mathbb R^n, 0<r\leq R} \ \left(r^{-n}\int_{P_r(x,r^2)}|\nabla u|^2+|\nabla v|^2\right)
^\frac12
\sup_{x\in\mathbb R^n,0<r\leq R} \left(r^{-n}\int_{P_r(x,r^2)}|\nabla(u-v)|^2\right)^\frac12\\
&\leq& C\epsilon_2 \left[\sup_{0<t\leq R^2} \ \sqrt{t}\|\nabla (u-v)(t)\|_{L^\infty(\mathbb R^n)}
+\sup_{x\in\mathbb R^n, 0<r\leq R}
\left(r^{-n}\int_{P_r(x,r^2)}|\nabla u|^2+|\nabla v|^2\right)
^\frac12\right]\\
&\leq& C\epsilon_2|||u-v|||_{X_{R^2}}.
\end{eqnarray*}
\begin{eqnarray*}
&&II\\
&=&\sup_{0<t\leq R^2} \ t\||\nabla v|^2|u-v|(t)|\|_{L^\infty(\mathbb R^n)}
+\sup_{x\in\mathbb R^n, 0<r\leq R} \ r^{-n}\int_{P_r(x,r^2)}|\nabla v|^2|u-v|\\
&\leq& \left[\sup_{0<t\leq R^2}\sqrt t\|\nabla v(t)\|_{L^\infty(\mathbb R^n)}
+\sup_{x\in\mathbb R^n, 0<r\leq R} \ (r^{-n}\int_{P_r(x,r^2)}|\nabla v|^2)^\frac12\right]^2
\sup_{0<t\leq R^2}\|(u-v)(t)\|_{L^\infty(\mathbb R^n)}\\
&\leq& C \|v\|_{X_{R^2}}^2\sup_{0<t\leq R^2}\|(u-v)(t)\|_{L^\infty(\mathbb R^n)}
\leq C_4\epsilon_2^2 |||u-v|||_{X_{R^2}},
\end{eqnarray*}
where we have used Lemma 3.2 in the last step. Putting these two estimates together yields
$$|||\mathbf Tu-\mathbf Tv|||_{X_{R^2}}
\leq C(1+\epsilon_2)\epsilon_2|||u-v|||_{X_{R^2}}
\leq\theta_0 |||u-v|||_{X_{R^2}}$$
for some $\theta_0=\theta_0(\epsilon_2)\in (0,1)$, provided $\epsilon_2>0$ is sufficiently small.
\qed

\noindent{\bf Proof of Theorem 1.3}. It follows from Lemma 3.3, 3.4, and the fixed point theorem
that there exists $\epsilon_0=\epsilon_0(n,N)>0$ such that if $[u_0]_{{\rm{BMO}}_R(\mathbb R^n)}\le\epsilon_0$ for some $R>0$,
then there exists a unique $u\in X_{R^2}$ such that
$$u=\tilde{u}_0+\mathbb{S}(A(u)(\nabla u,\nabla u))\ \ {\rm{on}}\ \mathbb R^n\times [0,R^2],$$
or equivalently 
$$u_t-\Delta u = A(u)(\nabla u, \nabla u) \ {\rm{on}}\ \mathbb R^n\times (0,R^2);
\ u\big|_{t=0}=u_0.$$ 

Now we need to show $u(\mathbb R^n\times [0,R^2])\subset N$. First, observe that
Lemma 2.1 implies that for $0<t\le \frac{R^2}{K^2}$, 
\begin{eqnarray*}
{\rm{dist}}(u,N)&\leq &{\rm{dist}}(\tilde{u}_0,N)+\|u-\tilde{u}_0\|_{L^\infty(\mathbb R^n\times [0, \frac{R^2}{K^2}))}\\
&\leq& \delta +K^n\left[u_0\right]_{{\rm{BMO}}_R(\mathbb R^n)}+\epsilon_0\\
&\leq &\delta+(1+K^n)\epsilon_0\le\delta_N,
\end{eqnarray*}
provide $\delta\le \frac{\delta_N}2$ and $\epsilon_0\le \frac{\delta_N}{2(1+K^n)}$.  
This yields $u(\mathbb R^n\times [0,\frac{R^2}{K^2}])
\subset N_{\delta_N}$. 
This and the definition of $A(\cdot)(\cdot,\cdot)$ imply  
$$A(u)(\nabla u,\nabla u)=-\nabla^2\Pi(u)(\nabla u,\nabla u) \ {\rm{on}}\ \mathbb R^n\times 
[0,\frac{R^2}{K^2}].$$
Set $Q(y)=y-\Pi(y)$ for $y\in N_{\delta_N}$,
and $\rho(u)=\frac12|Q(u)|^2$. Then direct calculations imply that for any $y\in N_{\delta_N}$,
$$\nabla Q(y)(v)=({\rm{Id}}-\nabla\Pi(y))(v), \ \forall v\in \mathbb R^l,$$
and
$$\nabla^2 Q(y)(v,w)=-\nabla^2\Pi(y)(v,w), \ \forall v, w\in\mathbb R^l.$$
Hence we have
\begin{eqnarray}
&&(\partial_t-\Delta) \rho(u)\nonumber\\
&=&\langle Q(u), \nabla Q(u)(\partial_t u-\Delta u)-\nabla^2 Q(u)(\nabla u,\nabla u)\rangle
   -|\nabla(Q(u))|^2\nonumber\\
&=&\langle Q(u), -\nabla Q(u)(\nabla^2\Pi(u)(\nabla u,\nabla u))-\nabla^2 Q(u)(\nabla u,\nabla u)\rangle
   -|\nabla(Q(u))|^2\nonumber\\
&=&\langle Q(u), \nabla \Pi(u)(\nabla^2\Pi(u)(\nabla u,\nabla u))\rangle -|\nabla (Q(u))|^2\nonumber\\
&=&-|\nabla(Q(u))|^2\le 0, \label{subharm}
\end{eqnarray}
where we have used the fact that $Q(u)\perp T_{\Pi(u)}N$ and $
\nabla \Pi(u)(\nabla^2\Pi(\nabla u,\nabla u))\in T_{\Pi(u)}N$ in the last step. 

Since $\rho(u)|_{t=0}=0$, the maximum principle for (\ref{subharm}) implies $\rho(u)\equiv 0$
on $\mathbb R^n\times [0,\frac{R^2}{K^2}])$.  One can repeat the same argument to show
that $u(\mathbb R^n\times [\frac{R^2}{K^2}, R^2])\subset N$. Thus
the proof of Theorem 1.3 is complete. \qed\\

\noindent{\bf Proof of Theorem 1.4}. It follows directly from Theorem 1.3 
with $R$ replaced by $+\infty$. \qed

\setcounter{section}{3} \setcounter{equation}{0}
\section{Proof of Theorem 1.5 and 1.6}

This section is devoted to the proof of Theorem 1.5 and 1.6 on local and global well-posedness of hydrodynamic flow of liquid crystals.

For $(u_0,d_0):\mathbb R^n\to\mathbb R^n\times S^2$, let $(\tilde{u}_0, \tilde{d}_0):
\mathbb R^n\times \mathbb R_+\to\mathbb R^n\times \mathbb R^3$ denote the
caloric extension of $(u_0,d_0)$. 

First, we recall the Carleson's characterization of $u_0\in {\rm{BMO}}_R^{-1}(\mathbb R^n)$
for $R>0$, due to Koch-Tataru \cite{KT}, which asserts that the following is equivalent
\beq{}\label{bmo-1}
[u_0]_{{\rm{BMO}_R^{-1}(\mathbb R^n)}}
\approx \sup_{x\in\mathbb R^n, 0<r\le R} (r^{-n}\int_{P_r(x,r^2)}|\tilde{u}_0|^2)^\frac12.
\eeq 

Notice that since $\tilde{u}_0$ solves the heat equation on
$\mathbb R^n$, the Harnack estimate of heat equation implies that
\begin{equation}\label{carloric_est}
\sup_{0<t\le R^2} \sqrt{t}\|\tilde{u}_0\|_{L^\infty}
\lesssim \sup_{x\in\mathbb R^n, 0<r\le R}\left(r^{-n}\int_{P_r(x,r^2)}|{\tilde u}_0|^2\right)^\frac12
\approx [u_0]_{\rm{BMO}^{-1}(\mathbb R^n)}.
\end{equation}

In particular, $u_0\in {\rm{BMO}}_R^{-1}(\mathbb R^n)$ implies that
$\tilde{u}_0\in Z_{R^2}$ and
\beq{}\label{Z-est}
\|\tilde{u}_0\|_{Z_{R^2}}\lesssim \|u_0\|_{{\rm{BMO}}_R^{-1}(\mathbb R^n)}.
\eeq

Let $\mathbb P: L^2(\mathbb R^n)\to \mathbb PL^2(\mathbb R^n)$ denote the Leray projection
operator. Then (\ref{llf1})-(\ref{llf2}) and $u|_{t=0}=u_0$
is equivalent to
\begin{equation}\label{llf5}
u(t)={\mathbb T}_1[u,d](t):={\tilde u}_0(t)-\mathbb{V}[u\otimes u+\nabla d\otimes \nabla d](t),
\end{equation}
where the operator $\mathbb V$ is defined by
\begin{equation}\label{v_op}
\mathbb{V}f(t)=\int_0^t e^{-(t-s)\Delta}\mathbb P\nabla\cdot f(s)\,ds,
\ \forall f:\mathbb R^n\times \mathbb R_+\to\mathbb R^n.
\end{equation}

The following estimate on the operator $\mathbb{V}$
has been proved by Koch-Tataru ([KT] Lemma 3.2).
\begin{Lemma} For $0<T\le+\infty$, if $f=(f_1,\cdots,f_n)\in Y_T$, then
\begin{equation} \label{XZ_bound}
\|Vf\|_{Z_{T}}\leq C\|f\|_{Y_T}
\end{equation}
for some constant $C=C(n)>0$.
\end{Lemma}

Observe that (\ref{llf3}) and $d|_{t=0}=d_0$ is equivalent to
\begin{equation}\label{llf6}
d(t)={\mathbb T}_2[u,d](t):={\tilde d}_0(t)+\mathbb{S}[-\nabla^2\Pi_{S^2}(d)(\nabla d,\nabla d)-u\cdot\nabla d](t),
\end{equation}
where $\mathbb S$ is the operator defined by (3.1), and $\Pi_{S^2}\in C^\infty(\mathbb R^3,\mathbb R^3)$
has the property 
$$\Pi_{S^2}(d)=\frac{d}{|d|}: S^2_{\frac12}\equiv\{y\in\mathbb R^3: \frac12\le|y|\le \frac32\}
\to S^2.$$ 

Let $(u_0,d_0)\in {\rm{BMO}}^{-1}_R(\mathbb R^n)\times {\rm{BMO}}_R(\mathbb R^n)$
for some $R>0$.  Then $({\tilde u}_0, {\tilde d}_0) \in 
Z_{R^2}\times X_{R^2}$. For $\epsilon>0$, we define 
the ball $\mathbb B_\epsilon([{\tilde u}_0, {\tilde d}_0])$ in $Z_{R^2}\times X_{R^2}$
with center $({\tilde u}_0, {\tilde d}_0)$ and radius $\epsilon$ by
$$
\mathbb B_\epsilon([{\tilde u}_0, {\tilde d}_0])
=\left\{(u,d)\in Z_{R^2}\times X_{R^2}: 
\left\|u-{\tilde u}_0\right\|_{Z_{R^2}}
+|||d-{\tilde d}_0|||_{X_{R^2}}\le\epsilon\right\}.
$$
Define the mapping operator $T$ on $Z_{R^2}\times X_{R^2}$ by
$${\mathbb T}[u,d]=({\mathbb T}_1[u,d], {\mathbb T}_2[u,d]).$$

Analogous to Lemma 3.2 and 3.3, we have the following two Lemmas.
\begin{Lemma} There exists $\epsilon_1>0$ such that if
$$\|u_0\|_{{\rm{BMO}}_R^{-1}(\mathbb R^n)}+\left[d_0\right]_{{\rm{BMO}}_R(\mathbb R^n)}\le\epsilon_1$$
then $\mathbb T$ maps  $\mathbb B_{\epsilon_1}([{\tilde u}_0, {\tilde d}_0])$
to $\mathbb B_{\epsilon_1}([{\tilde u}_0, {\tilde d}_0])$.
\end{Lemma}
\noindent{\it Proof}. For $(u,d)\in \mathbb B_{\epsilon_1}([{\tilde u}_0, {\tilde d}_0])$,
we have that $\|d\|_{L^\infty(\mathbb R^n\times [0,R^2])}\le C$ and
$$\mathbb T[u,d]-({\tilde u}_0, {\tilde d}_0)
=\left(-\mathbb{V}[u\otimes u+\nabla d\otimes \nabla d],\ 
\mathbb{S}[-\nabla^2 \Pi_{S^2}(d)(\nabla d,\nabla d)-u\cdot\nabla d]\right).
$$
Therefore, applying Lemma 3.1 and Lemma 4.1, we have
\begin{eqnarray*}
 &&\|\mathbb T_1[u,d]-{\tilde u}_0\|_{Z_{R^2}}
+|||\mathbb T_2[u,d]- {\tilde d}_0|||_{X_{R^2}}\\
&\lesssim&
\|u\otimes u+\nabla d\otimes \nabla d\|_{Y_{R^2}}+
\|\nabla^2 \Pi_{S^2}(d)(\nabla d,\nabla d)-u\cdot\nabla d\|_{Y_{R^2}}\\
&\lesssim&
\left(\left\|u\right\|_{Z_{R^2}}+\left\|d\right\|_{X_{R^2}}\right)^2\\
&\lesssim&
\left(\left\|u-{\tilde u_0}\right\|_{Z_{R^2}}+\left\|d-{\tilde d}_0\right\|_{X_{R^2}}
+\left\|{\tilde u}_0\right\|_{Z_{R^2}}+\left\|{\tilde d}_0\right\|_{X_{R^2}}\right)^2\\
&\leq& C\epsilon_1^2\le \epsilon_1
\end{eqnarray*}
provided $\epsilon_1>0$ is chosen to be sufficiently small, where we have
used the estimate 
$$\left\|{\tilde u}_0\right\|_{Z_{R^2}}+\left\|{\tilde d}_0\right\|_{X_{R^2}}
\lesssim \|u_0\|_{{\rm{BMO}}_R^{-1}(\mathbb R^n)}
+\left[d_0\right]_{{\rm{BMO}}_R(\mathbb R^n)}
$$
in the last step.
\qed

\begin{Lemma}
There exist $0<\epsilon_2\leq\epsilon_1$ and $\theta_0\in (0,1)$ such that if
$$\|u_0\|_{{\rm{BMO}}_R^{-1}(\mathbb R^n)}+\left[d_0\right]_{{\rm{BMO}}_R(\mathbb R^n)}\le\epsilon_2$$
then $\mathbb T: \mathbb B_{\epsilon_2}([{\tilde u}_0, {\tilde d}_0])
\to \mathbb B_{\epsilon_2}([{\tilde u}_0, {\tilde d}_0])$ is  $\theta_0$-contractive, i.e.
$$\|\mathbb T_1[u_1,d_1]-\mathbb T_1[u_2,d_2]\|_{Z_{R^2}}
+|||\mathbb T_2[u_1,d_1]-\mathbb T_2[u_2,d_2]|||_{X_{R^2}}
\le\theta_0(\|u_1-u_2\|_{Z_{R^2}}+|||d_1-d_2|||_{X_{R^2}})$$
for any $ (u_1,d_1) \ (u_2,d_2)\in B_{\epsilon_2}([\tilde{u}_0, \tilde{d}_0])$.
\end{Lemma}
\noindent{\it Proof}. For any $(u_1,d_1) \ (u_2,d_2)\in
B_{\epsilon_2}([\tilde{u}_0, \tilde{d}_0])$,
we have
\begin{eqnarray*}
&&|\mathbb T_1[u_1,d_1]-\mathbb T_1[u_2,d_2]|\\
&=& |\mathbb{V}[u_1\otimes u_1+\nabla d_1\otimes \nabla d_1-u_2\otimes u_2-\nabla d_2\otimes \nabla d_2]|\\
&\lesssim&
\mathbb{V}((|u_1|+|u_2|)|u_1-u_2|+(|\nabla d_1|+|\nabla d_2|)|\nabla (d_1-d_2)|),
\end{eqnarray*}
and
\begin{eqnarray*}
&&|\mathbb T_2[u_1,d_1]-\mathbb T_1[u_2,d_2]|\\
&=&|\mathbb{S}[-\nabla^2 \Pi_{S^2}(d_1)(\nabla d_1,\nabla d_1)-u_1\cdot\nabla d_1
+\nabla^2 \Pi_{S^2}(d_2)(\nabla d_2,\nabla d_2)+u_2\cdot\nabla d_2]|\\
&\lesssim&
\mathbb{S}((|\nabla d_1|+|\nabla d_2|+|u_1|)|\nabla(d_1-d_2)|+|\nabla d_2|^2|d_1-d_2|
+|u_1-u_2||\nabla d_2|).
\end{eqnarray*}
Thus Lemma 3.1 and Lemm 4.1 imply
\begin{eqnarray*}
&&\| \mathbb T_1[u_1,d_1]-\mathbb T_1[u_2,d_2]\|_{Z_{R^2}}
+||| \mathbb T_2[u_1,d_1]-\mathbb T_2[u_2,d_2]|||_{X_{R^2}}\\
&\lesssim&
\|(|u_1|+|u_2|)|u_1-u_2|+(|\nabla d_1|+|\nabla d_2|)|\nabla (d_1-d_2)|\|_{Y_{R^2}}\\
&+&
\|(|\nabla d_1|+|\nabla d_2|+|u_1|)|\nabla(d_1-d_2)|+|\nabla d_2|^2|d_1-d_2|
+|u_1-u_2||\nabla d_2|\|_{Y_{R^2}}\\
&\leq&
C\epsilon_2 \left[\|u_1-u_2\|_{Z_{R^2}}+|||d_1-d_2|||_{X_{R^2}}\right]\\
&\leq& \theta_0\left[\|u_1-u_2\|_{Z_{R^2}}+|||d_1-d_2|||_{X_{R^2}}\right]
\end{eqnarray*}
for some $\theta_0\in (0,1)$, provided $\epsilon_2>0$ is chosen to be sufficiently small,
where we have used 
$$\|u_i\|_{Z_{R^2}}+\|d_i\|_{X_{R^2}}\le C\epsilon_2, \ i=1,2$$
in the last steps.
This completes the proof.
\qed\\

\noindent{\bf Proof of Theorem 1.5}. 
It follows directly from Lemma 4.2, Lemma 4.3, and the fixed point theory
that there exists $\epsilon_0>0$ such that if
$$\|u_0\|_{{\rm{BMO}}^{-1}_R(\mathbb R^n)}+[d_0]_{{\rm{BMO}}_R(\mathbb R^n)}
\le\epsilon_0,$$
then there exists $(u,d)\in Z_{R^2}\times X_{R^2}$ such that 
(\ref{llf1}), (\ref{llf2}), (\ref{llf4}), and (\ref{llf3}) replaced
by
\beq{}\label{llf7} 
d_t+u\cdot\nabla d-\Delta d=-\nabla^2\Pi_{S^2}(d)(\nabla d,\nabla d)
\eeq
hold. To complete the proof, we need to show $d(\mathbb R^n\times [0,R^2])\subset S^2$.
This step is  similar to  the proof of Theorem 1.3. First, Lemma 2.1 implies
that for $t\le \frac{R^2}{K^2}$, 
$${\rm{dist}}(d, S^2)\le \epsilon_0+\delta+K^2 [d_0]_{{\rm{BMO}}_R(\mathbb R^n)}
\le(1+K^n)\epsilon_0+\delta
\le\frac12,$$
provided $\delta\le\frac14$ and $\epsilon_0 \le\frac{1}{4(1+K^n)}$.  
Thus $d(\mathbb R^n\times [0,\frac{R^2}{K^2}])\subset S^2_{\frac12}$.
Now consider the function $\rho(d)=\frac12 |d-\Pi_{S^2}(d)|^2$. Then the same calculation
as in the proof of Theorem 1.3 gives
$$(\rho(d))_t+u\cdot\nabla (\rho(d))
-\Delta(\rho(d))=-|\nabla (d-\Pi_{S^2}(d))|^2\le 0.$$
Since $\rho(d)\big|_{t=0}=0$, the maximum principle implies $\rho(d)\equiv 0$
on $\mathbb R^n\times [0\frac{R^2}{K^2}]$ and $d(\mathbb R^n\times [0,
\frac{R^2}{K^2}])\subset S^2$.  Repeating the same argument
can imply $d(\mathbb R^n\times [\frac{R^2}{K^2},R^2))\subset S^2$.
The proof is complete.\qed\\

\noindent{\bf Proof of Theorem 1.6}.  It follows directly from Theorem 1.5 with $R$
replaced by $R=+\infty$.
\qed

\bigskip
\noindent{\bf Acknowledgements}. This work is partially supported by NSF grant 0601182.
The author is grateful to Professor Sverak for bringing the problem of well-posedness 
for the heat flow of harmonic maps to my attentions and offering several discussions, 
and his interests in this work.  The work is carried out while the author is visiting IMA 
as a New Directions Research Professorship. The author is indebited to  IMA for providing
both the financial support and the excellent research environment.

\end{document}